\newcommand \nid {\noindent}

\def \QED{\hfill \rule [0.25ex] {1ex}{1ex} \newline}

\newcommand \vsp[1]{\vspace{#1\baselineskip}}

\newenvironment {proof}{\par\nid{\bf Proof}\/:\ }{\vsp{0.25}\QED}

\def \benum {\begin {enumerate} }
\def \eenum {\end {enumerate} }
\def \bdesc {\begin {description} }
\def \edesc {\end {description} }
\def \bcent {\begin {center} }
\def \ecent {\end {center} }
\def \beqn {\begin {eqnarray} }
\def \eeqn {\end {eqnarray} }
\def \bitem {\begin {itemize} }
\def \eitem {\end {itemize} }

\documentclass[12pt]{article} 
\usepackage{latexsym} 

\newtheorem{claim}{Claim}[section]
\newtheorem{lemma}{Lemma}[section]
\newtheorem{theorem}{Theorem}[section]
\newtheorem{remark}{Remark}[section]
\newtheorem{fact}{Fact}[section]
\newtheorem{corr}{Corollary}[section]
\newtheorem{defn}{Definition}[section]
\def \bthm {\begin {theorem}}
\def \ethm {\end {theorem}}
\def \bprf {\begin {proof}}
\def \eprf {\end {proof}}
\def \blem {\begin {lemma}}
\def \elem {\end{lemma}}
\def \bclm {\begin{claim}}
\def \eclm {\end{claim}}
\def \brem {\begin{remark}}
\def \erem {\end{remark}}
\def \bcorr {\begin{corr}}
\def \ecorr {\end{corr}}
\def \bfact {\begin{fact}}
\def \efact {\end{fact}}
\def \bdefn {\begin{defn}}
\def \edefn {\end{defn}}
\def \bdesc {\begin{description}}
\def \edesc {\end{description}}

\def \ck1  {{\cal C}_{k}} 
\def \c2   {{\cal C}_{2}}

\begin{document} 

\title{Cartesian product and acyclic edge colouring}  
\author{Rahul Muthu ~~and~~  C.R. Subramanian \\  
The Institute of Mathematical Sciences, \\ 
Taramani, Chennai - 600 113, India. \\  
{\em email}:~\{rahulm,crs\}@imsc.res.in 
} 
\date{} 
\maketitle{} 

\begin{abstract}
  The {\em acyclic chromatic index}, denoted by $a'(G)$, of a graph 
  $G$ is the minimum number of colours used in any proper edge colouring
  of $G$ such that the union of any two colour classes does not   contain a cycle, that is, forms a forest. 
  We show that $a'(G \Box H) \leq a'(G)+a'(H)$ 
  for any two graphs $G$ and $H$ such that $\max\{a'(G),a'(H)\} > 1$. Here, 
  $G \Box H$ denotes the cartesian product of $G$ and $H$. 
  This extends a recent result of \cite{cocoon06} where 
  tight and constructive bounds on $a'(G)$ were obtained for a 
  class of grid-like graphs which can be expressed as the cartesian product 
  of a number of paths and cycles. 
  % efficiently. This is obtained as a corollary to an earlier theorem
  % relating the acyclic chromatic index of the {\em cartesian product}
  % of two graphs to the acyclic chromatic index of the first graph and
  % the maximum degree of the second graph. The limitation of the result
  % is that the second graph is required to be either an edge a path or
  % a cycle. Here we relax this restriction and a get a result where the
  % acyclic chromatic index of the resultant graph is stated in terms of
  % the acyclic chromatic indices of both constituent graphs.
\end{abstract} 

\section{Introduction}
All graphs we consider are simple and finite. Throughout the paper we
use $\Delta=\Delta(G)$ to denote the maximum degree of a graph $G$.
A colouring of the edges of a graph is {\em proper} if no pair
of incident edges receive the same colour. A proper colouring
$\mathcal{C}$ of the edges of a graph $G$ is {\em acyclic} if there is
no two-coloured (bichromatic) cycle in $G$ with respect to $\mathcal{C}$. 
In other words, the subgraph induced by the union of any two colour classes in
$\mathcal{C}$ is a forest. The minimum number of colours required to
edge-colour a graph $G$ acyclically is termed the {\em acyclic chromatic
index} of $G$ and is denoted by $a'(G)$. The notion of acyclic colouring was
introduced by Gr\"unbaum in \cite{grunb}.
The acyclic chromatic index and its vertex analogue are closely
related to other parameters like {\em oriented chromatic number} and 
{\em star chromatic number} of a graph $G$ both of which have many 
practical applications \cite{orientchrjgt97,starchrdm01}.

Determining $a'(G)$ is a hard problem both from a theoretical and 
from an algorithmic point of view. Even for the simple and highly structured 
class of complete graphs $(K_n)$, the value of $a'(G)$ is still not 
determined exactly. 

However, using probabilistic arguments, some loose upper bounds have
been obtained. For example, see

% $(i)$ ~\cite{aloncolouring} for a bound of $a'(G) \leq 64\Delta$ .

$(i)$ \cite{furtheraspect} for a bound of $a'(G) \leq 16\Delta$ which 
improves a previous bound of $a'(G) \leq 64\Delta$ due to 
\cite{aloncolouring}. 

$(ii)$ \cite{graco} for a bound of $a'(G) \leq 4.52\Delta$ for graphs $G$ with girth (the length of the shortest cycle) at least 220.

It has been conjectured  \cite{alonacyclic} that 
$a'(G) \leq \Delta+2$ for any $G$ and this has been shown to be true 
for some special classes of graphs. However, the presently known bounds 
are far from the conjectured bound of $\Delta+2$. It is still open 
whether the conjecture is true or if there are counterexamples.    

Some tight upper bounds have also been obtained for some special 
classes of graphs. For example, see

$(iii)$ \cite{burns} for a bound of $a'(G) \leq 5$ for 3-regular 
 graphs,   

$(iv)$ \cite{alonacyclic} for a bound of $a'(G) \leq \Delta+2$ for graphs 
with girth at least $c\Delta (\log \Delta)$, 

$(v)$ \cite{wornes} for a bound of $a'(G) \leq d+1$ for random  
$d$-regular ($d$ fixed) graphs.

Some constructive bounds which lead to an actual acyclic edge colouring have also been obtained. See 

$(vi)$ \cite{subcubic} for a constructive bound of 
$a'(G)\le 5$ for graphs with
$\Delta\le 3$,

$(vii)$ \cite{greedy} for a constructive bound of
$a'(G) \leq 6\Delta\log{\Delta}$ for any arbitrary graph,

$(viii)$ \cite{cocoon06} for a constructive bound of $\Delta+1$ for grid-like graphs. 

$(ix)$ \cite{aaim07} for a constructive bound of $\Delta+1$ for outerplanar graphs. 

% It has also been shown (see \cite{alonzaks}) that determining 
% if $a'(G) \leq 3$ is NP-complete for an arbitrary $G$. 

In this paper, we look at the cartesian product (defined in 
Section 2), denoted $G \Box H$, of two arbitrary graphs $G$ and $H$ and show that the acyclic chromatic index of the product is at most the sum of acyclic chromatic indices of $G$ and $H$, provided at least one of these values exceeds 1. This is an extension of the work in \cite{cocoon06}, where it was shown that 
$a'(G \Box H) \le a'(G)+\Delta(H)$,
whenever $H$ is a path or a cycle and $a'(G) \ge 2$, $a'(G) \ge 3$,
respectively. While the bound we give here is slightly weaker, we
remove the restriction that $H$ is a path or a cycle.
 
 Section $2$ contains definitions and our main result. Section $3$
contains some concluding remarks. In the following subsection, in order to motivate the reader, we  present a brief exposure to previous work on graph invariants in the context of graph products. 

\subsection{Previous work on graph invariants and cartesian product} 
Since it is well known (see \cite{imrichbook} for details) that any 
connected graph can be expressed in a unique way (upto isomorphism) as the cartesian product of the smaller so-called "prime" graphs, researchers have studied how various invariants of a graph can be expressed in terms of those of its factors. Specifically, it has been shown that 

 $(a)$ $\gamma(G \Box H) \leq \min\{\gamma(G) |V(G)|, \gamma(H) |V(H)|\}$ 
 by Vizing \cite{vizcarprod63}. Here, $\gamma(G)$ denotes the minimum size of a dominating set of $G$.  
 % It was also conjectured by Vizing that $\gamma(G \Box H) \geq 
 % \gamma(G) \gamma(H)$ and it is still open.  
 
 $(b)$ Vizing \cite{vizcarprod63} also studied the independence 
 number $\alpha(G)$ (the maximum size of an independent set of $G$) in the context of cartesian product of graphs and showed that 
 \begin{eqnarray}
\alpha(G \Box H) & \geq & \alpha(G) \alpha(H) + \min\{|V(G)|-\alpha(G), 
|V(H)|-\alpha(H)\} \nonumber \\ 
\alpha(G \Box H) & \leq & \min\{\alpha(G) |V(H)|, \alpha(H) |V(G)|\} \nonumber 
\end{eqnarray}

 $(c)$ Let $\chi(G)$ denote the chromatic number of $G$, that is, 
 the minimum number of colours required to properly colour the vertices of $G$. It was first noticed by Sabidussi \cite{sabidussi57} and 
 can also be easily verified that $\chi(G \Box H) = \max\{\chi(G), 
 \chi(H)\}$.  
 
  $(d)$ The Hadwiger number of a graph $G$, denoted by $\eta(G)$, is the largest integer $l$ such that $G$ has a $K_l$ minor where $K_l$ is the complete graph on $l $ vertices. Recently, Chandran and Raju \cite{sunilraju} have obtained results relating $\eta(G)$ and the cartesian product operation. In particular, it is proved in \cite{sunilraju} that $\eta(G \Box H) \ge (g-\sqrt{h})(\sqrt{h}-2)/4$ where 
  $g=\eta(G) \geq h=\eta(H)$.

\section{Definitions and Results}

For a comprehensive introduction and survey of results on various graph products, the reader is advised to refer to the book authored by 
Imrich and Klavzar \cite{imrichbook}. 
 
\bdefn Given two graphs $G_1=(V_1,E_1)$ and $G_2=(V_2,E_2)$, their
{\em cartesian product}, denoted by $G_1  \Box  G_2$, is defined as the graph
$G=(V,E)$ where $V=V_1 \times V_2$ and $E$ contains the edge joining
$(u_1,u_2)$ and $(v_1,v_2)$ if and only if {\em either} $u_1=v_1$ and
$\{u_2,v_2\}$ is an edge in $E_2$ {\em or} $u_2=v_2$ and 
$\{u_1,v_1\}$ is an edge in $E_1$.  
\edefn  

Note that $G=G_1 \Box G_2$ can be thought of as
being obtained as follows.  Take $|V_2|$ isomorphic copies of $G_1$
and label them with vertices from $V_2$. For each edge $\{u,v\}$ in
$E_2$, introduce a perfect matching between $G_u$ and $G_v$ which
joins each vertex in $V(G_u)$ with its isomorphic image in $V(G_v)$.
Equivalently, one can also think of this as being obtained by taking $|V_1|$
isomorphic copies of $G_2$ and introducing a perfect matching between
corresponding copies of $G_2$ for each edge in $E_1$.
The following facts are easy to verify. 

\bfact
The cartesian product $G_1  \Box  G_2$ is commutative 
in the sense that $G_1 \Box G_2$ is isomorphic to $G_2 \Box G_1$. Similarly, 
this operation is also associative. Hence the product 
$G_1 \Box G_2 \Box  \ldots  \Box G_k$ is well-defined for each $k$. For each 
$G$ and $k \geq 1$, we define $G^k$ as follows : $G^1=G$ and
$G^k = G^{k-1} \Box G$ for $k>1$.  
\efact  

\bfact
If $G=G_1 \Box G_2 \Box  \ldots  \Box G_k$, then $G=(V,E)$ where $V$ is the 
set of all $k$-tuples of the form $(u_1, \ldots ,u_k)$ with each 
$u_i \in V(G_i)$ and the edge joining 
$(u_1, \ldots ,u_k)$ and  $(v_1, \ldots ,v_k)$ is in $E$ if and only if 
for some $i$, $1 \leq i \leq k$, $(i)$ $u_j = v_j$ for all $j \not= i$ 
and $(ii)$ the edge $\{u_i,v_i\}$ is in $E(G_i)$.  
\efact 

\bfact \label{boxconfact} 
$G_1 \Box G_2$ is connected if and only if both $G_1$ and $G_2$ are 
connected.  
\efact 

\brem
It is known  % shown by Sabidussi 
% \cite{sabidussi60} and Vizing \cite{vizing63} 
(see \cite{imrichbook} for further details and references) 
that any connected graph $G$ can be expressed as a product $G=G_1 \Box \cdots \Box 
G_k$ of prime factors $G_i$. Here, a graph is said to be {\em prime} 
with respect to the $\Box$ operation if it is non-trivial and if  
it is not isomorphic to the product of two non-trivial graphs. A non-trivial graph is one having at least two vertices. Also, this 
factorisation is {\em unique} except for a re-ordering of the factors 
and is known as the {\em Unique Prime Factorisation (UPF)} of the graph. It is also known that the UPF of a graph can be computed in time polynomial in the 
size of $G$. 
\erem 

We now formally present our main result which relates acyclic chromatic index to the cartesian product of graphs. Without loss of generality, we can assume that the product graph is connected. In view of Fact \ref{boxconfact}, this implies that it suffices  
to consider only connected graphs as factors. 
Also, if $H$ is trivial (that is, 
$H$ is a graph on just one vertex), then $G \Box H$ is isomorphic to $G$ for 
any $G$. Hence, we focus only on connected non-trivial graphs.   
We will often use the following easy-to-verify fact about acyclic edge colourings. 

\bfact 
For any $\Delta > 1$, let $G$ be any $\Delta$-regular graph. Then, $a'(G) \geq \Delta+1$. 
\efact 

\bthm \label{bpaclthm}  
Let $G=(V_G,E_G)$ and $H=(V_H,E_H)$ be two connected non-trivial graphs 
such that $\max\{a'(G),a'(H)\} > 1$. Then, 
$$a'(G \Box H) \leq a'(G)+a'(H).$$ 
\ethm 
\noindent{\bf Note :} If $G$ and $H$ are both connected and 
non-trivial with $a'(G)=a'(H)=1$, then each of $G$ and $H$ is a $K_2$. 
In that case, $G \Box H=C_4$ where $C_4$ is a cycle on 4 vertices. 
Only in this case, we have $a'(G \Box H)=3$ whereas $a'(G)+a'(H)=2$.  
\bprf 
Let $a'(G)=\eta$ and $a'(H)=\beta$. Since $\Box$ is commutative,  without loss of generality, assume that $\eta \geq \beta$. 
Let $\Delta$ denote the maximum degree of $H$. Set $d$ to be $\Delta+1$ if $H$ is either a complete graph on $\Delta+1$ vertices or an 
odd cycle (in which case $\Delta=2$). Otherwise, set $d$ to be 
$\Delta$. % either 
% $\Delta+1$ or $\Delta$ depending on whether $H=K_{\Delta+1}$ or not. 
In any case, $H$ can be properly {\em vertex} 
coloured using colours from the set $[d]=\{0, \ldots d-1\}$. 
% We will handle the case when $H=K_{d+1}$ separately later.  

We know that $\beta=a'(H) \geq \Delta$ always. 
If $H=K_{\Delta+1}$, then (since $H$ is $\Delta$-regular) 
$a'(H) \geq \Delta+1$ (except when $H=K_2$). In both cases, 
$\eta \geq \beta \geq d$. If $H=K_2$, then $d=\Delta+1=2$ and 
$\eta \geq 2$ by assumption. In any case, we have $\eta \geq d$.  
Let $X_G : E_G \rightarrow [\eta]=\{0, \ldots ,\eta-1\}$ and 
$X_H : E_H \rightarrow [\beta']=\{0', \ldots ,(\beta-1)'\}$ be 
two acyclic edge colourings of $G$ and $H$ respectively using 
disjoint sets of colours. 

Each edge in $G \Box H$ is either $(i)$ an edge joining 
$(u_1,v)$ and $(u_2,v)$ for some $e=\{u_1,u_2\} \in E_G$ and $v \in V_H$ 
or $(ii)$ an edge joining    
$(u,v_1)$ and $(u,v_2)$ for some $f=\{v_1,v_2\} \in E_H$ and $u \in V_G$. 
We denote the former edges by $e_v$ (where $e \in E_G, v \in V_H$) 
and the latter edges by $f_u$ (where $f \in E_H, u \in V_G$).  
Note that each edge of $G \Box H$ lies either in some isomorphic copy $H_u$ 
of $H$ or in some isomorphic copy $G_v$ of $G$.  

For each $i \in \{0, \ldots ,d-1\}$, let $\sigma_i : [\eta]  
\rightarrow [\eta]$ be a bijection defined by 
$$ \sigma_i(j) = (j+i) \mbox{ mod } \eta, \;\;\; \forall j \in [\eta].$$

Since $\eta \geq d$, we notice that the bijections 
$\sigma_i (i \in [d])$ are mutually non-fixing, that is, for all 
$0 \leq i,k \leq d-1, \; i \not= k$ and for each $j \in [\eta]$,   
$\; \sigma_i(j) \not= \sigma_k(j)$.   

Let $Y_H : V_H \rightarrow \{0, \ldots ,d-1\}$ be a proper vertex 
colouring of $V_H$. We define a colouring of the edges of $G \Box H$  
based on the colourings $X_G,X_H$ and $Y_H$ as follows. 

For each edge in $E$ of the form $f_u$, where $f \in E_H$ and 
$u \in V_G$, we colour $f_u$ using the colour $X_H(f)$. 
Now consider any arbitrary edge of the form $e_v$ where $e \in E_G$ 
and $v \in V_H$. Let $i = Y_H(v)$ be the colour used by $Y_H$ on $v$.  
Colour $e_v$ using the colour $\sigma_i(X_G(e))$.  

In other words, edges $f_u$ in each isomorphic copy $H_u$ is coloured 
the same way as $f$ in $H$ is coloured by $X_H$. But edges $e_v$ in each 
isomorphic copy $G_v$ is coloured essentially (ignoring the labels of 
colours) the same way as $G$ is coloured but the colour labels are 
rotated by mutually non-fixing permutations. The permutation that 
is used for a $G_v$ is decided by the vertex colour assigned to $v$ 
by $Y_H$. As a result, for each edge $f=\{v_1,v_2\} \in E_H$ and for each 
edge $e=\{u_1,u_2\} \in E_G$, $e_{v_1}$ and $e_{v_2}$ get different colours 
but always from $[\eta]$.

Let $X : E(G \Box H) \rightarrow \{0, \ldots ,\eta-1\} \cup 
\{0', \ldots ,(\beta-1)'\}$ be the colouring defined just above. We 
will show that $X$ is proper and acyclic. 
\bclm
X is proper. 
\eclm  
\bprf
Consider any vertex $(u,v)$. The set of edges in $G \Box H$ which 
are incident on $(u,v)$ can be partitioned into two subsets 
$A_u = \{f_u : v \in f \in E_H\}$ and   
$A_v = \{e_v : u \in e \in E_G\}$. Since edges in these two sets 
are coloured using colours from disjoint sets, namely from 
$[\eta]$ and $[\beta']$, there is no conflict 
between these two sets. Now, let us focus on edges in $A_u$. 
Since $f_u$'s are coloured in the same way as $f$'s are coloured 
in $H$, there is no conflict among edges in $A_u$. Similarly, 
the edges $e_v$'s in $A_v$ are coloured essentially in the same way (except for a rotation of the colour labels) as $e$'s are coloured in  $G$ and hence coloured with distinct colours, there 
is no conflict among members of $A_v$ also. Hence $X$ is proper.  
\eprf 
It is only left to prove acyclicity of $X$. We prove this by contradiction. 
Suppose there is a bichromatic (with respect to $X$) cycle $C$ in $G$. 
First, we note that
\bclm
$C$ cannot lie entirely within any isomorphic copy $G_v$ or $H_u$ of 
$G$ or $H$ respectively.  
\eclm 
\bprf
Note that $X$ restricted to $H_u$ (or $G_v$) is basically either 
$X_H$ (or $X_G$ except for renaming of the colours). Hence if 
$C$ lies within such an isomorphic copy, it implies that either 
$X_H$ or $X_G$ has a bichromatic cycle, which is a contradiction. 
\eprf 
By the above claim, it follows that $C$ should visit vertices in 
at least two different copies $G_{v}$ and $G_{v'}$. But different copies 
are only joined by edges of type $f_u$ for some $f \in E_H$ and 
$u \in V_G$. Thus, it follows that $C$ has at least one edge each 
of the two types $e_v$ ($e \in E_G, \; v \in V_H$) and   
$f_u$ ($f \in E_H, \; u \in V_G$) which are coloured with respectively, 
say, $a \in [\eta]$ and $b \in [\beta']$. 

\bclm 
Let $(u_1,v_1)$ be some arbitrary vertex in $C$. Let $(u_1,v_2)$ 
for some $v_2 \in V_H$ be the other end point of the unique $b$-coloured 
edge in $C$  incident at $(u_1,v_1)$. 
$C$ lies entirely within $G_{v_1}$ and $G_{v_2}$. 
\eclm  
\bprf
The proof is by induction on the distance $l$ in $C$ from $(u_1,v_1)$ 
{\em along the direction} specified by the edge $\{v_1,v_2\}_{u_1}$. 
For $l=0$, it is clearly true. Suppose it is true for vertices whose 
above-defined distance is at most $l'$. Let $(u_{l'},v_{l'})$ be 
the vertex at distance $l'$. By inductive hypothesis, $v_{l'}$ is 
either $v_1$ or $v_2$. Let $c \in \{a,b\}$ be the colour of the 
edge joining $(u_{l'},v_{l'})$ and $(u_{l'+1},v_{l'+1})$. If 
$c=a$, then $v_{l'+1}=v_{l'}$ and hence the hypothesis is clearly 
true for $l=l'+1$. If $c=b$ (hence $u_{l'+1}=u_{l'}$) and if 
$v_{l'}=v_1$, then $v_{l'+1} = v_2$. 
This follows from $(i)$ the $b$-coloured edge incident at the 
copy of $u_1$ in $G_{v_1}$ joins it to the copy of $u_1$ in $G_{v_2}$ 
and hence $(ii)$ all edges of the perfect matching joining isomorphic 
copies of vertices in $G_{v_1}$ and $G_{v_2}$ are coloured with 
$b$. In particular, the $b$-coloured edge incident at $(u_{l'},v_1)$ 
joins it to $(u_{l'},v_2)$. 
Similarly, one can argue that if $c=b$ and if $v_{l'}=v_2$, then 
$v_{l'+1} = v_1$. In any case, $v_{l'+1} \in \{v_1,v_2\}$, thereby 
proving that $C$ lies entirely within $G_{v_1}$ and $G_{v_2}$.  
\eprf 
Since the edges in $G_{v_1}$ and $G_{v_2}$ are coloured without using 
colour $b$ and since every alternate edge of $C$ is coloured with $b$, 
we see that $b$ is used an even number of times in $C$. This implies 
$|C| = 0 \; (\mbox{mod} \; 4)$. 
Thus, $C$ looks like $$C = \langle \; (u_1,v_1), (u_1,v_2), (u_2,v_2), 
(u_2,v_1), \ldots, (u_{2k-1},v_2), (u_{2k},v_2), 
(u_{2k},v_1), (u_1,v_1) \; \rangle.$$   
For each of the $a$-coloured edges in $G_{v_2}$ joining 
$(u_{2l-1},v_2)$ and $(u_{2l},v_2)$, its isomorphic copy in $G_{v_1}$ 
joins $(u_{2l-1},v_1)$ and $(u_{2l},v_1)$ and is coloured with 
the colour $c=\sigma_i(\sigma^{-1}_{j}(a)) \not= a$ where 
$i=Y_H(v_1)$ and $j=Y_H(v_2)$. These isomorphic copies in $G_{v_1}$ 
of $a$-coloured edges of $C$ in $G_{v_2}$ together with the $a$-coloured 
edges of $C$ in $G_{v_1}$ constitute the following $\{a,c\}$-coloured bichromatic cycle in $G_{v_1}$ : 
$$D = \langle \; (u_1,v_1), (u_2,v_1), (u_3,v_1), \ldots, 
(u_{2k},v_1), (u_1,v_1) \; \rangle.$$ 
This is a contradiction to the fact that $X$ restricted to $G_{v_1}$ 
is acyclic. This shows that $X$ admits no bichromatic cycle and hence 
$X$ is proper and acyclic. Since $X$ uses only colours from 
$[\eta] \cup [\beta']$, we get $a'(G \Box H) \leq a'(G)+a'(H)$. 
\eprf 
\brem
Note that the above proof is constructive in the following sense : given two acyclic edge colourings $X_G$ and $X_H$ of $G$ and $H$ respectively, one can construct an acyclic edge colouring of $G \Box H$ in time polynomial in the size of $E(G \Box H)$. 
\erem 
\bcorr
Let $G_1, \ldots ,G_k$ be $k$ connected non-trivial graphs such that for 
each $i$, $1 \leq i \leq k$, $a'(G_i)=\Delta(G_i)$ and 
$\max\{a'(G_1), \ldots ,a'(G_k)\} >1$. Then,  
$$a'(G_1 \Box  \ldots  \Box G_k) = \Delta(G_1 \Box  \ldots  \Box G_k).$$   
\ecorr 
\bprf
Follows from $(i)$ $a'(G) \geq \Delta(G)$ for any $G$, 
$(ii)$ $\Delta(G_1 \Box  \ldots  \Box G_k) = \Delta(G_1)+ \ldots +\Delta(G_k)$, 
$(iii)$ Theorem \ref{bpaclthm}. 
\eprf 
\bcorr
Let $G$ be a connected non-trivial graph such that $a'(G)=\Delta(G) > 1$.  
Then, for each $d \geq 1$,  
$$a'(G^d) = d\Delta(G).$$   
\ecorr
The following result first obtained in \cite{cocoon06} now follows 
as a corollary of Theorem 2.1. 
\bcorr (\cite{cocoon06}) 
Let $G=K_2^d=K_2 \Box  \ldots  \Box K_2$ be the $d$-dimensional hypercube for 
some $d \geq 1$. Then, 
$$a'(K_2) = 1 \;\;\; \mbox{ and } \;\;\; a'(K_2^d)=d+1 \;\; \mbox{ for } 
\;\; d>1.$$  
\ecorr   
\bprf
Suppose $d>1$. Since $G=K_2^d$ is $d$-regular, by Fact 2.4, we need at least $d+1$ colours in any acyclic edge colouring of $K_2^d$ and hence 
$a'(G) \geq d+1$. Also, $a'(K_2^2)=a'(C_4)=3$. Starting with 
$G=K_2^2$ and applying Theorem \ref{bpaclthm} repeatedly by 
setting $H=K_2$ each time, we get $a'(K_2^d) \leq a'(K_2^2)+(d-2) 
\leq d+1$. Combining both the lower and upper bounds, we get the 
result.   
\eprf

\section{Conclusions}
It can be easily observed that $a'(G)\ge \Delta+1$, for all regular
graphs with $\Delta\ge2$. It is conjectured in \cite{alonacyclic} that
$a'(G)\le\Delta+2$ for all graphs. If we take the cartesian product of
$t$ graphs each of whose $a'(G)$ value is known and each of which is regular with $\Delta \ge 2$, then the
bound we would get (which is not better than $\Delta+t$) by applying the above
result is very weak, {\em assuming} the conjecture is true. It would
be interesting to make a statement like $a'(G \Box H)\le a'(G)+\Delta(H)$
for a wider class of graphs $H$, like the results obtained in \cite{cocoon06} for grid-like graphs. At present, investigations are being carried out in this direction. 
% \bthm \label{bpaclconjthm}  
% Let $G=(V_G,E_G)$ and $H=(V_H,E_H)$ be two connected non-trivial graphs 
% such that $(i)$ $a'(G) \leq \Delta(G)+2$ and 
% $a'(H) \leq \Delta(H)+2$; $(ii)$ $\max\{a'(G),a'(H)\} > 1$. Then, 
% $$a'(G \Box H) \leq \Delta(G \Box H)+2.$$ 
% \ethm 

\end{document}